\newif\ifpictures
\picturestrue

\documentclass[12pt]{amsart}
\usepackage{amssymb}
\usepackage[dvips]{graphicx}

\ifpictures
\usepackage{psfrag}
\fi

\numberwithin{equation}{section}

\newtheorem{thm}{Theorem}
\newtheorem{prop}[thm]{Proposition}
\newtheorem{lemma}[thm]{Lemma}
\newtheorem{cor}[thm]{Corollary}

\theoremstyle{definition}

\newtheorem{example}[thm]{Example}
\newtheorem{remark1}[thm]{Remark}
\newtheorem{openproblem1}[thm]{Open problem}

\newenvironment{rem}{\begin{remark1}\rm}{\end{remark1}}

\newenvironment{openproblem}{\begin{openproblem1}\rm}{\end{openproblem1}}

\numberwithin{thm}{section}

\newcounter{FNC}[page]
\def\newfootnote#1{{\addtocounter{FNC}{2}$^\fnsymbol{FNC}$%
     \let\thefootnote\relax\footnotetext{$^\fnsymbol{FNC}$#1}}}

\newcommand{\C}{\mathbb{C}}
\newcommand{\N}{\mathbb{N}}
\newcommand{\Q}{\mathbb{Q}}
\newcommand{\R}{\mathbb{R}}

\newcommand{\classP}{\mathcal{P}}
\newcommand{\NP}{\mathcal{N}\mathcal{P}}

\DeclareMathOperator{\conv}{conv}

\DeclareMathOperator{\trop}{trop}

\DeclareMathOperator{\ord}{ord}

\def\problem#1#2#3{ 

\medskip

\noindent {\bf Problem {#1}:} \medskip

\noindent
\begin{tabular}{@{}lp{12cm}}
  {\bf Instance:} & {#2} \\
  {\bf Question:} & {#3} \\
\end{tabular}

\bigskip

}

\title[Computations in tropical geometry]{On the frontiers of polynomial
computations \\ in tropical geometry}

\author{Thorsten Theobald}
\address{Technische Universit\"at Berlin, Stra{ss}e des 17.~Juni 136,
  D-10623 Berlin, Germany}
\email{theobald@math.tu-berlin.de}

\thanks{Part of this work was done while the author was a Feodor Lynen fellow of
the Alexander von Humboldt Foundation at Yale University.}
\begin{document}

\begin{abstract}
We study some basic algorithmic problems concerning the intersection of
tropical hypersurfaces in general dimension: deciding whether 
this intersection is non\-emp\-ty, whether it is a tropical variety,
and whether it is connected, as well as counting the number
of connected components. We characterize the borderline between
tractable and hard computations by proving $\mathcal{NP}$-hardness
and \#$\mathcal{P}$-hardness
results under various strong restrictions of the
input data, as well as providing polynomial time algorithms for
various other restrictions. 
\end{abstract}

\maketitle




\section{Introduction}

Geometry over the tropical semiring 
$(R, \oplus, \odot) := (\R,\min,+)$
has received much attention in the last years (see the 
surveys \cite{mikhalkin-survey-2004,rst-2004,speyer-sturmfels-survey-2004}
and the references therein)
with applications in counting curves \cite{mikhalkin-2005},
studying phylogenetic trees \cite{pachter-sturmfels-2004},
and the analysis of amoebas of algebraic 
varieties \cite{mikhalkin-survey-2004}.
From the viewpoint of polynomial equations, the modern
birth of tropical geometry originates
in the book \cite{sturmfels-b2002} which pinpoints the
central role of tropical geometry as a link between
algebraic geometry, symbolic computation, and discrete geometry,
thus providing computationally-accessible methods
for studying algebraic-geometric problems. Indeed, one
of the early roots of the developments in
tropical geometry can be seen in the polyhedral homotopy methods
by Huber and Sturmfels~\cite{huber-sturmfels-95}, providing
a state-of-the-art technique for numerically solving systems of polynomial
equations based on a deformation to a (discrete) tropical problem.

Some major algorithmic results in tropical geometry
are based on Gr\"obner basis computations
and thus may become intractable already for small dimensions
\cite{bjsst-05}.
On the positive side, there also exist some algorithmic problems (such as computing
the tropical determinant) which
can be efficiently solved using techniques from linear
programming, polyhedral computation and combinatorial optimization 
(see, e.g., \cite{joswig-2003,rst-2004}). For many tropical problems,
their computational complexity has not been clarified yet.

In this paper, we make a first step towards systematically studying
the frontiers of polynomial time computations
in tropical
geometry. For this, we consider three natural algorithmic problems
concerning the intersection of tropical hypersurfaces, so-called
tropical prevarieties. The algorithmic problems are to decide whether
this intersection $P$ is nonempty (\textsc{Tropical Intersection}), 
whether $P$ is a tropical variety (\textsc{Tropical Consistency}), 
and whether $P$ is connected (\textsc{Tropical Connectivity}).

Our results refer to the the standard Turing machine model, 
and we mainly aim at characterizing the borderline between
tractable (in the sense of polynomial time solvable) and
hard (in the sense of $\NP$-hard) computations.
Our main results can roughly be stated as follows.
If the number of hypersurfaces is part of the 
input then the three problems become $\NP$-hard or co-$\NP$-hard, 
and this hardness
persists even under various restrictions to the input data.
As a particular example, already for quadratic input polynomials it is 
co-$\NP$-hard to decide whether a tropical 
prevariety is a tropical variety.
Hence, efficient algorithms cannot be expected in this setting.
We contrast these hardness results by polynomial time algorithms 
for a fixed number of tropical hypersurfaces.
For a precise statement of the results see 
Theorems~\ref{th:hardness1}--\ref{th:sharpp}.

The paper is structured as follows.
In Section~\ref{se:prelim} we introduce the relevant notation from
tropical geometry and computational complexity.
In Section~\ref{se:mainresults} we formally state our main results,
and Section~\ref{se:proofs} contains the proofs of these
theorems.
We close the paper with a short
discussion of related computational aspects on amoebas.

\section{preliminaries\label{se:prelim}}

\subsection{Tropical geometry}

One of the original motivations for tropical varieties was a combinatorial
approach to certain problems from enumerative geometry suggested by Kontsevich, 
and that program has been realized by Mikhalkin \cite{mikhalkin-2005}.
Tropical varieties are also related to the observation that algebraic
varieties have very simple behavior at infinity when plotted on ``log paper''
\cite{bergman-71, viro-2001}. While these roots come from algebraic geometry
and valuation theory, tropical varieties are profitably approached via 
polyhedral combinatorics.

Tropical hypersurfaces can be defined in a combinatorial and
in an algebraic way (for general background we refer 
to \cite{mikhalkin-survey-2004}, \cite{rst-2004},
\cite[Chapter 9]{sturmfels-b2002}).
For the combinatorial definition, let
$(\R, \oplus, \odot)$ denote the \emph{tropical semiring}, where
\[
  x \oplus y \ = \ \min \{x,y\} \quad \text{ and } \quad
  x \odot y \ = \ x + y \, .
\]
Sometimes the underlying set $\R$ of real numbers is augmented by $\infty$.

A \emph{tropical monomial} is an expression of the form
$c \odot x^{\alpha} = c \odot x_1^{\alpha_1} \odot \cdots \odot x_n^{\alpha_n}$
where the powers of the variables are computed tropically as well
(e.g., $x_1^3 = x_1 \odot x_1 \odot x_1)$. This tropical monomial
represents the classical linear function
\[
  \R^n \to \R \, , \quad 
  (x_1, \ldots, x_n) \mapsto \alpha_1 x_1 + \dots + \alpha_n x_n + c \, .
\]

A \emph{tropical polynomial} is a finite tropical sum of tropical monomials
and thus represents the (pointwise) minimum function of linear functions.
At each given point $x \in \R^n$ the minimum is either attained at a single
linear function or at more than one of the linear functions
(``at least twice'').
The \emph{tropical hypersurface} $\mathcal{T}(f)$ of a tropical polynomial $f$
is defined by
\[
  \begin{array}{rcl}
  \mathcal{T}(f) \ = \ \{ x \in \R^n & : & \text{ the minimum of the tropical monomials of $f$} \\
  & & \text{ is attained at least twice at $x$} \} 
  \, .
  \end{array}
\]

Rather than simply intersecting tropical hypersurfaces, the definition of
tropical varieties of arbitrary codimension involves a valuation
theoretic construction (Section~\ref{se:geometrytropical} explains this
subtlety.)
Let $K = \overline{\C(t)}$ denote the algebraically closed field of 
Puiseux series, i.e., series of the form
\[
   p(t) \ = \ c_1 t^{q_1} \,+\, c_2 t^{q_2}  \,+ \, c_3 t^{q_3} \,+ \,
\cdots
\]
with $c_i \in \C \setminus \{0\}$ and rational
 $q_1 < q_2 < \cdots $ with common denominator
 (see, e.g., \cite{bpr-2003}).
 The \emph{order} $\ord p(t)$ is the exponent of the lowest-order term
 of $p(t)$.
 The order of an $n$-tuple of Puiseux series
 is the $n$-tuple of their orders. This gives a map
 \begin{equation}
 \label{ordmap}
  {\rm ord} \,\,: \,\, (K^*)^n \,\, \rightarrow \,\,
 \Q^n \,\,\, \subset  \,\,\, \R^n \, ,
 \end{equation}
 where $K^* = K \setminus \{0\}$.
 
We are extending $\mathcal{T}$ to allow also ideals in 
the polynomial ring $K[x_1, \ldots, x_n]$ as argument.
Let $I$ be an ideal in $K[x_1, \ldots, x_n]$, and consider its
affine variety $V(I) \subset  (K^*)^n$
over $K$.  The image of
$V(I)$ under the map (\ref{ordmap}) is a subset
of $\Q^n$. The \emph{tropical variety} $\mathcal{T}(I)$
is defined as the topological closure of this image,
$\mathcal{T}(I) = \overline{\ord V(I)}$.
It is well-known that for principal ideals 
$I = \langle g \rangle$ the two definitions 
of tropical varieties coincide (see 
\cite{kapranov-2000} or, e.g., \cite[Lemma 3.2]{rst-2004}):

\begin{prop} If $f$ is a tropical polynomial in $x_1, \ldots, x_n$
then there exists
a polynomial $g \in K[x_1, \ldots, x_n]$ such that
$\mathcal{T}(f) = \mathcal{T}(\langle g \rangle)$, and vice versa.
\end{prop}

For a polynomial $f = \sum_{\alpha \in \mathcal{A}} c_{\alpha}(t) x^{\alpha}
\in K[x_1, \ldots, x_n]$ with a finite support set $\mathcal{A} \subset \N_0^n$
and $c_{\alpha}(t) \neq 0$ for all $\alpha \in \mathcal{A}$, the 
\emph{tropicalization} of $f$ is defined by
\[
  \trop f \ = \ \bigoplus_{\alpha \in \mathcal{A}} \ord(c_{\alpha}(t)) \odot
  x^{\alpha} \, ,
\]
where $\bigoplus$ denotes a tropical summation.
Whenever there is no possibility of confusion we also 
write $\cdot$ instead of $\odot$.

For every tropical variety $\mathcal{T}(I)$ there exists a finite
subset $\mathcal{B} \subset I$ such that
$\mathcal{T}(I) = \bigcap_{f \in \mathcal{B}} \mathcal{T}(f)$.
(However, we remark that Corollary 2.3 in \cite{speyer-sturmfels-2004},
which claims that any universal Gr\"obner basis of $I$ satisfies
this condition, is not correct. See \cite{rst-2004}.)

\ifpictures
\begin{figure}[ht]
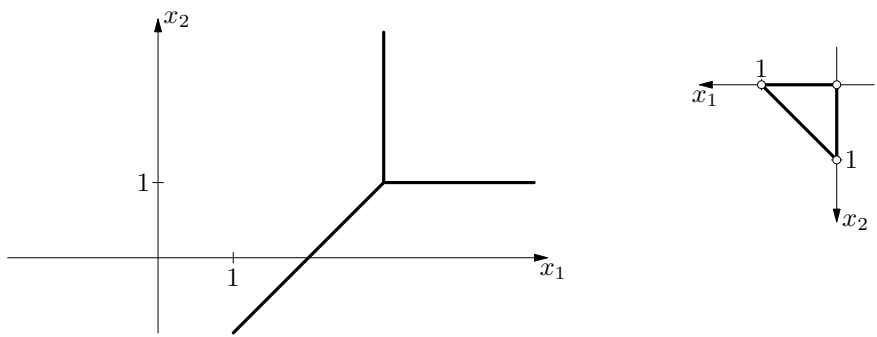

\[
  \begin{array}{c@{\qquad \qquad}c}
    \includegraphics[scale=1]{pictures/pictropcomp.6} \\ [-4cm] 
    & \includegraphics[scale=1]{pictures/pictropcomp.7} \\ [1cm]
  \end{array}
\]

\caption{The tropical variety of a linear polynomial $f$ in two
  variables and the Newton polygon of $f$.}
\label{fi:linearcurve1}
\end{figure}
\fi

\subsection{The geometry of tropical hypersurfaces\label{se:geometrytropical}}
Let $\mathcal{A} \subset \N_0^n$ be finite and
$f(x_1, \ldots, x_n) = \bigoplus_{\alpha \in \mathcal{A}} c_{\alpha} \cdot x^{\alpha}$
be a tropical polynomial
with $c_{\alpha} \in \R$ for all $\alpha \in \mathcal{A}$. Then
$\mathcal{T}(f)$ is a polyhedral complex in $\R^n$ 
which is geometrically dual to
the following regular subdivision of the Newton 
polytope $\text{New}(f)$ of $f$.
Let $\hat{P}$ be the convex hull
$\conv \{ (\alpha, c_{\alpha}) \in \R^{n+1} \, : \, 
  \alpha \in \mathcal{A} \}$.
Then the lower faces of $\hat{P}$ project bijectively onto 
$\conv \mathcal{A}$ under deletion of the last coordinate,
thus defining a subdivision of $\mathcal{A}$. Such subdivisions are
called \emph{regular} or \emph{coherent} (see, e.g., \cite{lee-91}).
We say that a tropical polynomial is \emph{of degree at most $d$} 
if every term has (total) degree at most~$d$.
See Figure~\ref{fi:linearcurve1} for an example of a tropical line
(i.e., the tropical variety of a linear polynomial in two variables) and 
Figure~\ref{fi:tropcubic1} for an example of a tropical cubic curve,
as well as their dual subdivisions (whose coordinate axes are directed 
to the left and to the bottom to visualize the duality).

\ifpictures
\begin{figure}[ht]
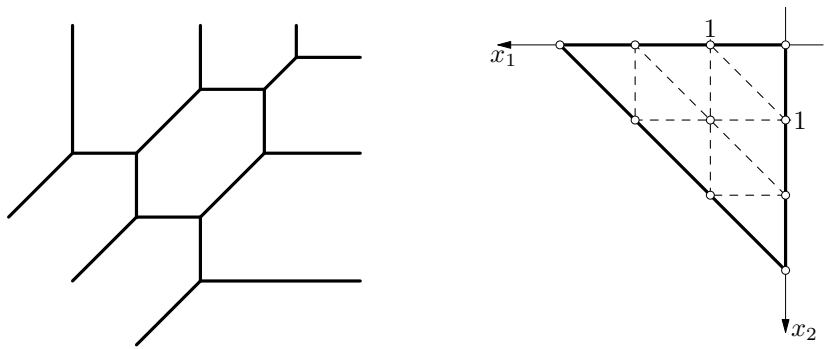

\[
  \begin{array}{c@{\qquad \qquad}c}
    \includegraphics[scale=1]{pictures/pictropcomp.8} \\ [-4.7cm] 
    & \includegraphics[scale=1]{pictures/pictropcomp.9} \\ [0.2cm]
  \end{array}
\]

\caption{An example of a tropical cubic curve $\mathcal{T}(f)$
  and the dual subdivision of the Newton polygon of $f$.}
\label{fi:tropcubic1}
\end{figure}
\fi

Following the notation in \cite{rst-2004},
a \emph{tropical prevariety} is the intersection of tropical
hypersurfaces.
If $f_1, \ldots, f_m$ are linear polynomials then the tropical
prevariety $P = \bigcap_{i=1}^m \mathcal{T}(f_i)$ is called \emph{linear}.
If additionally $P$ is a tropical variety, then it
is called a \emph{linear tropical variety}.
In dimension~2, a linear tropical variety is either
a translate of the left-hand set
in Figure~\ref{fi:linearcurve1}, a classical line (in the $x_1$-, $x_2$-, or
the main diagonal direction), a single point, or
the empty set. A tropical prevariety in $\R^2$ can also be a 
one-sided infinite ray.
Understanding the geometry and combinatorics
of tropical prevarieties or varieties in general dimension
is still a widely open problem. 
Even for the case of linear tropical varieties, the maximum number of bounded 
$i$-dimensional faces of such polyhedral complexes is unknown.
The recent
\emph{f-vector conjecture} in \cite{speyer-2004}
conjectures that (in our affine setting) the number of bounded
$i$-dimensional faces of a $k$-dimensional linear tropical variety
in 
$\R^n$ is at most $\binom{n-2i+1}{k-i+1}\binom{n-i}{i-1}$
and that this bound is tight.

With respect to our investigations on the consistency problem, we remark
that there are linear tropical spaces of dimension $n-2$
which are not complete intersections, i.e., which are not the 
intersection of two tropical hypersurfaces 
(see \cite[Proposition 6.3]{speyer-sturmfels-2004}).

\subsection{Model of computation.}

Our model of computation is the binary Turing machine:
all relevant data are presented by certain rational numbers, and the size
of the input is defined as the length of the binary encoding of
the input data. A rational number is specified as the concatenation 
of the numerator $a$ and the denominator $b$, and we may assume without loss of
generality that $a$ and $b$ are relatively prime.
Polynomials of degree~$d$ are specified
by the binary encoding of all $\binom{n+d}{n}$ coefficients (even if a coefficient
is zero); this encoding is sometimes referred as the \emph{dense} encoding.
For general background on algorithms and complexity theory 
we refer to \cite{garey-johnson-b79,papadimitriou-steiglitz-82}, 
and in particular for
complexity aspects of geometric problems to \cite{gls-b93}.

In the realm of the complexity classes $\classP$ and $\NP$,
complexity theory usually deals with decision problems: those whose answers are
$\textsc{Yes}$ or $\textsc{No}$.
The class $\classP$ denotes the set of all decision problems which
can be solved in polynomial time in the input size.
The class $\NP$ (nondeterministic polynomial time) denotes the 
class of all problems such that every \textsc{Yes}-instance has a short
(i.e.\ polynomial-size) certificate that can be verified in polynomial
time. Recall that a problem is called \emph{co-$\mathcal{NP}$-hard} if its
complement is $\mathcal{NP}$-hard, where the complement of a problem
is defined by switching the answers $\textsc{Yes}$ and $\textsc{No}$ for all
inputs.

In this paper, we also deal with counting problems, which refer
to problems whose answer is a bit string encoding an integer.
A counting problem $\Pi$ is in
the class \#$\mathcal{P}$ if there is a decision problem $\Pi' \in \mathcal{NP}$
such that, for all inputs $I$, the output of $\Pi$ is exactly
the number of accepting solutions to $\Pi'$ on input $I$.
A counting problem $\Pi$ is \#$\mathcal{P}$-\emph{hard} if every problem in
\#$\classP$ can be reduced in polynomial time to $\Pi$, i.e.,
if for every problem $\Pi \in \#\mathcal{P}$ 
there is a polynomial-time computable function $f$ such that
for any input $I$ to $\Pi'$
\begin{enumerate}
\item $f(I)$ is a valid input to $\Pi$, 
\item the output of $\Pi'$ on input $I$ is exactly the output of $\Pi$ on input
  $f(I)$.
\end{enumerate}

\section{Statement of problems and main results\label{se:mainresults}}

We consider three basic problems on the intersection of tropical
hypersurfaces. Let $\Q[x_1, \ldots, x_n]^{\oplus}$ denote the set
of tropical polynomials with rational coefficients in $n$ variables.
Given $n$, $m$, $d_1, \ldots, d_m$ and a set of tropical polynomials 
$f_1, \ldots, f_m \in \Q[x_1, \ldots, x_n]^{\oplus}$ of degrees at
most $d_1, \ldots, d_m$, respectively, the first problem
asks whether the tropical prevariety $\bigcap_{i=1}^m \mathcal{T}(f_i)$
is nonempty. For the complexity results it is quite crucial
which information is part of the input of the problem. 
In particular, note that in the formal definitions of the
three problems the dimension and the number of hypersurfaces
is part of the input.

\problem{\rm \textsc{Tropical Intersection}}
  {$n$, $m$, $d_1, \ldots, d_m$, 
  polynomials $f_1, \ldots, f_m \in \Q[x_1, \ldots, x_n]^{\oplus}$
  of degrees at most $d_1, \ldots, d_m$.}
  {Decide whether there exists a point in 
  $\bigcap\limits_{i=1}^m \mathcal{T}(f_i)$.}

The next problem asks whether an intersection of tropical hypersurfaces
(i.e., a prevariety) is a tropical variety.

\problem{\rm \textsc{Tropical Consistency}}
  {$n$, $m$, $d_1, \ldots, d_m$,
  polynomials $f_1, \ldots, f_m \in \Q[x_1, \ldots, x_n]^{\oplus}$
  of degrees at most $d_1, \ldots, d_m$.}
  {Decide whether $\bigcap\limits_{i=1}^m \mathcal{T}(f_i)$ is a
   tropical variety.}

We also consider the variant $\textsc{Tropical $m$-Consistency}$ which asks
whether $\bigcap_{i=1}^m \mathcal{T}(f_i)$ is a tropical variety
of codimension $m$.
The third problem asks for topological connectivity of the set 
$\bigcap_{i=1}^m \mathcal{T}(f_i)$.

\problem{\rm \textsc{Tropical Connectivity}}
  {$n$, $m$, $d_1, \ldots, d_m$,
  polynomials $f_1, \ldots, f_m \in \Q[x_1, \ldots, x_n]^{\oplus}$ 
  of degrees at most $d_1, \ldots, d_m$ with
  $\bigcap\limits_{i=1}^m \mathcal{T}(f_i) \neq \emptyset$.}
  {Decide whether $\bigcap\limits_{i=1}^m \mathcal{T}(f_i)$ is connected.}

Besides these decision problems, we consider the counting problem
\textsc{\#Connected Components} whose input is the same one
as for \textsc{Tropical Intersection} and whose task is to determine
the number of connected components of 
$\bigcap_{i=1}^m \mathcal{T}(f_i)$.

Our main results can be stated as follows.

\begin{thm}
\label{th:hardness1}
The problem \textsc{Tropical Intersection} is $\NP$-complete,
and the problems \textsc{Tropical Consistency} and
\textsc{Tropical Connectivity} are co-$\NP$-hard.
For \textsc{Tropical Intersection} and \textsc{Tropical Connectivity}
these hardness results persist if the instances are restricted to those
where $\bigcap_{i=1}^m \mathcal{T}(f_i)$ is a 
tropical variety.

Moreover, for \textsc{Tropical Intersection} and \textsc{Tropical Consistency}
the hardness persists if all polynomials are restricted to
be of degree at most~2.
For \textsc{Tropical Connectivity},
the hardness persists if all polynomials are restricted to
be of degree at most~3.
\end{thm}

These hardness results are contrasted by the following positive algorithmic
results for restricted input classes.

\begin{thm}
\label{th:easyness1}
(i) If the number $m$ of tropical hypersurfaces is a fixed constant, then
\textsc{Tropical Intersection} can be solved in polynomial time.
\smallskip

\noindent
(ii) For fixed $m$ and if all input polynomials are restricted to be linear
   polynomials then the problem \textsc{Tropical $m$-Consistency}
   can be solved in polynomial time.
\smallskip

\noindent
(iii) If the number $m$ of tropical hypersurfaces is a fixed constant, then
\textsc{Tropical Connectivity} can be solved in polynomial time.
Moreover, 
any linear tropical prevariety is connected; hence, if all polynomials
are restricted to be linear polynomials, the output of
\textsc{Tropical Connectivity} is always \textsc{Yes}.
\end{thm}

Finally, we show \#$\classP$-hardness of counting the number of solutions.

\begin{thm}
\label{th:sharpp}
$\#$\textsc{Connected Components} is $\#\mathcal{P}$-hard. 
This statement persists
if all polynomials are restricted to be of degree at most 2.
\end{thm}

\begin{rem}
Obviously, \textsc{Tropical Intersection} (and similarly, 
\textsc{Tropical Connectivity} and \textsc{Connected Components})
can be solved (not necessarily efficiently)
by explicitly constructing the polyhedral complexes 
$\mathcal{T}(f_1), \ldots, \mathcal{T}(f_m)$ in $\R^n$
and intersecting them. 

Solving \textsc{Tropical Consistency} in a similar way
can be done based on a synthetic definition of the tropical varieties
under investigation.
For tropical hypersurfaces such a definition can be found 
in~\cite[Prop.~3.15]{mikhalkin-2005} and for tropical lines in
$\R^n$ in \cite[Example 3.8]{rst-2004}.
\end{rem}

Several questions remain open. In particular, 
the question of polynomial time solvability remains open 
for the following restrictions.

\begin{openproblem}
Can  
\textsc{Tropical Intersection} and \textsc{Tropical Consistency}
be solved in polynomial time if the input polynomials are restricted
to be linear?
Can \textsc{Tropical Connectivity} for quadratic polynomials be solved
in polynomial time?
\end{openproblem}

\section{Proofs of the results\label{se:proofs}}

\subsection{Linear tropical prevarieties\label{se:lineartropical}}

We begin with a statement on \textsc{Tropical Consistency and}
\textsc{Tropical $m$-Consistency}
for linear varieties.

\begin{lemma} \label{le:mconsistency}
Let all input polynomials $f_1, \ldots, f_m \in \Q[x_1, \ldots, x_n]^{\oplus}$
be restricted to be linear polynomials.
\begin{enumerate}
\item[(a)] If $m \le n$ then the output of \textsc{Tropical Consistency} is always {\sc Yes}.
\item[(b)] For a fixed constant $m$, the problems \textsc{Tropical $m$-Consistency}
can be solved in polynomial time.
\end{enumerate}
\end{lemma}

Before providing the proof, we recall and collect some statements
about linear tropical varieties.
Let $f_1, \ldots, f_m$ be linear tropical polynomials in $x_1, \ldots, x_n$.
If $m \le n$ and the tropical hyperplanes $\mathcal{T}(f_i)$ are
in general position then $P$ is a linear tropical variety
of dimension $n-m$. Moreover, $P$ always contains a well-defined
\emph{stable intersection} which is a linear tropical variety
of dimension $n-m$ (see \cite{rst-2004,speyer-2004}). In particular,
this implies that for $m \le n$ the answer to \textsc{Tropical Consistency} is
always {\sc Yes.}

For a matrix $A = (a_{ij}) \in (\R \cup \{\infty\})^{k \times k}$,
the tropical determinant is defined by
\begin{equation}
  \label{eq:det}
  \mathrm{det}_{\text{trop}}(A) \ = \ 
  \bigoplus_{\sigma \in S_k} (a_{1,\sigma_1} \odot \cdots \odot a_{k,\sigma_k})
  \ = \ \min_{\sigma \in S_k} (a_{1,\sigma_1} + \cdots + a_{k,\sigma_k}) \, ,
\end{equation}
where $S_k$ denotes the symmetric group on $\{1, \ldots, k\}$. It
was observed in \cite{rst-2004} that the computation of the tropical
determinant can be phrased as an assignment problem from combinatorial
optimization. Hence, using well-known 
algorithms (see \cite[Corollary 17.4b]{schrijver-2003}),
a tropical determinant can be computed in polynomial time.

A tropical $k \times k$-matrix is \emph{singular} if the minimum
in~\eqref{eq:det} is attained at least twice. In order to
decide in polynomial time whether a $k \times k$-matrix is singular,
first compute the tropical determinant of $A$.
Let $\sigma \in S_k$ be a permutation of $\{1, \ldots, k\}$ for which
the minimum in~\eqref{eq:det} is attained. For every 
$j \in \{1, \ldots, k\}$ let $A_j$ be the matrix which is obtained 
from $A$ by replacing the entry $(j,\sigma_j)$ by an arbitrary larger value.
Then $A$ is tropically singular if and only if 
$\det_{\text{trop}} A_j = \det_{\text{trop}} A$ for some
$j \in \{1, \ldots, k\}$. Hence, this can be decided in polynomial time.

\medskip

\noindent
\emph{Proof of Lemma~\ref{le:mconsistency}.} \smallskip
It just remains to prove b).
Let $f_i = \bigoplus_{j=1}^n a_{ij} \cdot x_j \oplus a_{i,n+1}$, 
$1 \le i \le m$, and $A = (a_{ij})\in \R^{m \times (n+1)}$ be the 
coefficient matrix of $f_1, \ldots, f_m$. Since for $m > n+1$ 
the answer of \textsc{Tropical $m$-Consistency} is always 
$\textsc{No}$, we can assume $m \le n+1$. For $m=n+1$, the problem
is equivalent to ask whether the tropical prevariety is empty,
which will be treated in Lemma~\ref{le:fixedm}.
For $m \le n$,
by Theorem~5.3 in \cite{rst-2004} the tropical prevariety
$\bigcap_{i=1}^m \mathcal{T}(f_i)$
is a linear tropical variety of codimension $m$ if and only if
none of the $m \times m$-submatrices of $A$ is tropically singular.
For each of the $m \times m$-submatrices of $A$
it can be checked in polynomial time (in the binary length of
the input data) whether it is singular.
Since for fixed $m$, the number $\binom{n}{m}$
of those submatrices is polynomial in $n$, the claim follows.
\hfill $\Box$

\subsection{Tropical intersection and tropical consistency\label{se:tropintersection}}

\begin{lemma} \label{le:tropintersection}
{\sc Tropical Intersection} is $\NP$-hard.
This statement persists if the instances are restricted to
those where $\bigcap_{i=1}^m \mathcal{T}(f_i)$ is a tropical variety.
Moreover,
this statement persists if all polynomials are restricted to be
of degree at most 2.
\end{lemma}

In order to prove $\NP$-hardness of {\sc Tropical Intersection}, 
we provide a polynomial time reduction from the
well-known $\NP$-complete 3-satisfiability (3-{\sc Sat})
problem \cite{garey-johnson-b79}.
Let $\wedge$ and $\vee$ denote the Boolean conjunction and disjunction,
respectively, and let
$\mathcal{C} = \mathcal{C}_1 \wedge \ldots \wedge \mathcal{C}_k$
denote an instance of 3-{\sc Sat} with clauses
$\mathcal{C}_1, \ldots, \mathcal{C}_k$ in the variables
$y_1, \ldots, y_n$.
Furthermore, let $\overline{y_i}$ denote the complement of a variable
$y_i$, and let the literals $y_i^1$ and $y_i^{0}$ be defined by
$y_i^1 = y_i$, $y_i^{0} = \overline{y_i}$. Let
the clause $\mathcal{C}_i$ be of the form
\begin{equation}
\label{eq:clause1}
  \mathcal{C}_i = y_{i_1}^{\tau_{i_1}} \vee y_{i_2}^{\tau_{i_2}} \vee 
  y_{i_3}^{\tau_{i_3}},
\end{equation}
where $\tau_{i_1}, \tau_{i_2}, \tau_{i_3} \in \{0,1\}$ and
$i_1, i_2, i_3 \in \{1, \ldots, n\}$ are pairwise different indices.

The reduction consists of two ingredients.
First we construct an intersection of suitable tropical hypersurfaces
$\bigcap_{i=1}^n \mathcal{T}(h_i)$ in $\R^n$
with $\bigcap_{i=1}^n \mathcal{T}(h_i) = \{0,1\}^n$
(see Figure~\ref{fi:structuralhypersurfaces}).
We call these hypersurfaces
``structural'' tropical hypersurfaces.

\ifpictures
\begin{figure}[ht]
\[
  \includegraphics[scale=1]{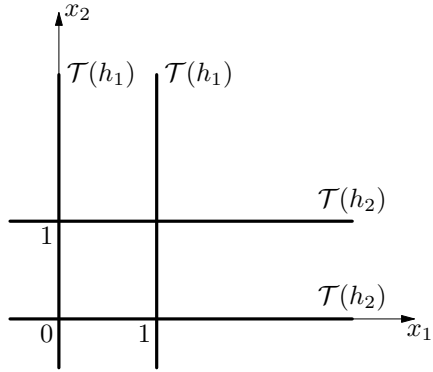}
\]

\caption{Structural hypersurfaces 
$\mathcal{T}(h_i)$, $1 \le i \le n$, for $n=2$.}
\label{fi:structuralhypersurfaces}
\end{figure}
\fi

In the second step, we relate satisfying assignments of a given
clause~\eqref{eq:clause1} to solutions
of some ``clause hypersurface''.
Let $s:\{ \text{\sc True}, \text{\sc False}\} \to \{0,1\}$ be defined
by $s(\text{\sc True}) = 1$ and
$s(\text{\sc False}) = 0$. 
We utilize the correspondence between a truth assignment
$a = (a_1, \ldots, a_n)^T \in \{\text{\sc True},\text{\sc False}\}^n$
to the variables $y_1, \ldots, y_n$
and the point $(s(a_1), \ldots, s(a_n))^T \in \{0,1\}^n$
of the tropical prevariety $\bigcap_{i=1}^n \mathcal{T}(h_i)$. 
To achieve this, we construct one or,
in some cases, several tropical hypersurfaces representing the clause.

In order to  construct the structural tropical hypersurfaces,
let $h_i' \in K[x_1, \ldots, x_n]$ be the polynomial
\[
  h_i'(x) \ = \ 
  (t^0 \cdot x_i + t^1) \cdot (t^0 \cdot x_i + t^{0})
  \ = \ t^0 \cdot x_i^2 + (t^0 + t^1) \cdot x_i + t^1
\]
over $K$, $1 \le i \le n$.
Since the tropical hypersurface of a product of polynomials is
the union of the tropical hypersurfaces of the factors, we have
$\mathcal{T}(h_i') = \left\{ x \in \R^n \, : \, x_i \in \{0,1\} \right\}$,
and $h_i'$ tropicalizes to
\begin{equation}
\label{eq:structuralclause1}
  h_i \ := \text{trop}(h_i') \ = \ 0 \cdot x_i^2 \oplus 0 \cdot x_i \oplus 1 \, .
\end{equation}
Hence, $\bigcap_{i=1}^n \mathcal{T}(h_i) \ = \ \{0,1\}^n$. 

Now we construct the quadratic polynomials which represent the
3-clauses.
In order to illustrate the
construction, and since this will be needed explicitly later on, we begin
with a 2-clause. 
Let $\mathcal{C}_i$ denote the 2-clause
$\mathcal{C}_i = y_{i_1}^{\tau_{i_1}} \vee y_{i_2}^{\tau_{i_2}}$. 
Let $f_i' = (x_{i_1} + t^{\tau_{i_1}})(x_{i_2} + t^{\tau_{i_2}})$, 
which tropicalizes to 
\begin{eqnarray*}
  f_i & = & (0 \cdot x_{i_1} \oplus \tau_{i_1}) \cdot 
                         (0 \cdot x_{i_2} \oplus \tau_{i_2}) \\
  & = & 0 \cdot x_{i_1} \cdot x_{i_2} \oplus \tau_{i_1} \cdot x_{i_1} 
        \oplus \tau_{i_2} \cdot x_{i_2}
        \oplus \tau_{i_1} \cdot \tau_{i_2} \, .
\end{eqnarray*}
Hence, $\mathcal{T}(f_i) = \{x \in \R^2 \, : \, x_{i_1} = \tau_{i_1} 
\text{ or } x_{i_2} = \tau_{i_2} \}$.
In particular, for any $x \in \{0,1\}^n$ we have
$x \in \mathcal{T}(f_i)$ if and only if
$x_{i_1} = \tau_{i_1}$ or $x_{i_2} = \tau_{i_2}$.
Figure~\ref{fi:twoclause} (a) shows the clause curve for the
case $n=2$, $\tau_{1} = \tau_{2} = 1$.

Now consider a 3-clause $\mathcal{C}_i$ of the form \eqref{eq:clause1}.
Here, the straightforward approach to consider
$\mathcal{T}(f_i) = \prod_{j=1}^3 (0 \cdot x_{i_j} \oplus \tau_{i_j})$
leads to cubic polynomials. In order to show hardness even for 
quadratic polynomials
we distinguish several cases corresponding to the number $p$ of 
positive literals in $\mathcal{C}_i$.

\medskip

\noindent
\emph{Case $p \in \{0,1\}$:} Here we can use the following
more general lemma.

\begin{lemma} \label{le:clausetotropical1}
Let $C(y_1, \ldots, y_k, z_1, \ldots, z_l)$ be the clause
in the variables $y_1, \ldots, y_k, z_1, \ldots, z_l$ defined by
\[
  C(y_1, \ldots,y_k, z_1, \ldots, z_l) \ = \ y_1 \vee \cdots \vee y_k \vee
          \overline{z_1} \vee \cdots \vee \overline{z_l} \, .
\]
Then for $(a_1, \ldots, a_k, b_1, \ldots, b_l) \in \{\textsc{True},
\textsc{False}\}^{k+l}$
we have $C(a_1, \ldots, a_k, b_1, \ldots, b_l) = \textsc{True}$ if
and only if
\begin{equation}
  \label{eq:clausesurface1}
  (s(a_1), \ldots, s(a_k), s(b_1), \ldots, s(b_l)) \in
  \mathcal{T}\left( \left\langle \prod_{i=1}^k \left(t^0 \cdot y_i + t^1 \right) 
         \cdot \left( \sum_{j=1}^l t^0 \cdot z_j + t^{0} \right) \right\rangle \right) \, .
\end{equation}
\end{lemma}

\begin{proof}
Let $C_i(y_i) = y_i$, $1 \le i \le k$. Then
for $a_i \in \{\textsc{True},\textsc{False}\}$, we have
$C_i(a_i) = \textsc{True}$ if and only if
$s(a_i) \in \mathcal{T}( \left\langle t^0 \cdot y_i + t^1 \right\rangle )$.
Let $C_{k+1}(z_1, \ldots, z_l) = \overline{z_1} \vee \cdots
  \vee \overline{z_l}$. Then 
for $(b_1, \ldots, b_l) \in \{\textsc{True},\textsc{False} \}^l$, 
we have $C_{k+1}(z_1, \ldots, z_l) = \textsc{True}$ if and only if
$s(b_1, \ldots, b_l) \in 
\mathcal{T}(\langle \sum_{j=1}^l t^0 \cdot z_j + t^{0} \rangle )$.
Considering the disjunction $C_1 \vee \dots \vee C_{k+1}$
proves the claim.
\end{proof}

For every clause $\mathcal{C}_i$ which contains
0 or 1 positive literals, we associate a tropical 
hypersurface $\mathcal{T}(f_i)$ as defined in~\eqref{eq:clausesurface1}.
Since $p \in \{0,1\}$ the degree of $f_i$ is at most 2.
 
In particular, for the case $p = 0$ 
and $i_1 = 1$, $i_2 = 2$, $i_3 = 3$, we have
$\mathcal{C}_i = \overline{y_1} \vee \overline{y_2} \vee \overline{y_3}$,
and the hypersurface in~\eqref{eq:clausesurface1} is 
$\mathcal{T}(\langle t^0 \cdot y_1 + t^0 \cdot y_2 + t^0 \cdot y_3 + t^0 \rangle)$,
which is the hypersurface given by the linear tropical polynomial
$0 \cdot y_1 \oplus 0 \cdot y_2 \oplus 0 \cdot y_3$.
Figure~\ref{fi:twoclause}(b) visualizes this situation for the
smaller-dimensional case of a 2-clause.

\ifpictures
\begin{figure}[ht]
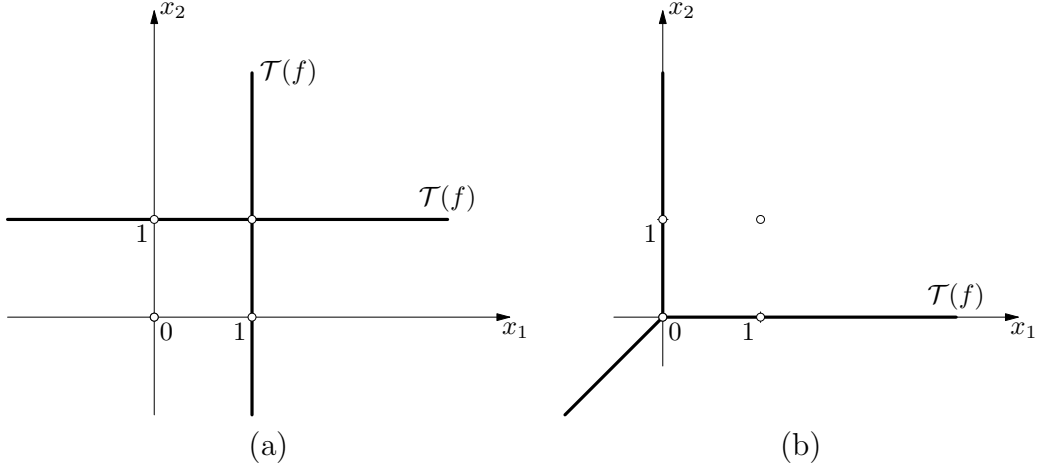

\[
  \begin{array}{c@{\quad}c}
    \includegraphics[scale=1]{pictures/pictropcomp.2} &
    \includegraphics[scale=1]{pictures/pictropcomp.3} \\
    \text{(a)} & \text{(b)}
  \end{array}
\]

\caption{Figure (a) shows a quadratic clause curve for the 2-clause
  $y_1 \vee y_2$ in $\R^2$. Figure (b) shows a linear clause curve
  for the 2-clause $\overline{y_1} \vee \overline{y_2}$ in $\R^2$.}
\label{fi:twoclause}
\end{figure}
\fi

For the case $p=1$ and the clause
$y_1 \vee \overline{y_2} \vee \overline{y_3}$, 
the hypersurface in~\eqref{eq:clausesurface1} is
\[
  \mathcal{T}\left( \left\langle (t^0 \cdot y_1 + t^1) \cdot 
  (t^0 \cdot y_2 + t^0 \cdot y_3 + t^0) \right\rangle \right) \, ,
\]
which is the hypersurface of the tropical polynomial
$0 \cdot y_1 \cdot y_2 \oplus 0 \cdot y_1 \cdot y_3 \oplus
0 \cdot y_1 \oplus 1 \cdot y_2 \oplus 1 \cdot y_3 \oplus 1$.

\medskip

\noindent
\emph{Case $p=2$:} By renumbering the variables, we can assume
$\mathcal{C}_i = y_{i_1} \vee y_{i_2} \vee \overline{y_{i_3}}$.
Let $f_i$ be the quadratic tropical polynomial defined by
\[
  f_i \ = \ 0 \cdot x_{i_1} \cdot x_{i_2} 
  \oplus 1 \cdot x_{i_1} 
  \oplus 1 \cdot x_{i_2}
  \oplus 0 \cdot x_{i_3}
  \oplus 0 \cdot x_{i_3}^2 \oplus 1 \, .
\]
Then for $(x_{i_1}, x_{i_2}, x_{i_3}) \in \{0,1\}^3$ we have
$f_i(x_{i_1},x_{i_2}, x_{i_3}) = 1$ if and only if 
$x_{i_1} = 1$ or $x_{i_2} = 1$ or $x_{i_3} = 0$.

\medskip

\noindent
\emph{Case $p=3$:} Let $\mathcal{C}_i = y_{i_1} \vee y_{i_2} \vee y_{i_3}$.
The following lemma (in the spirit of the nine associated points
theorem for complex cubic curves) 
states that it is \emph{not} possible to find 
a single polynomial $f_i$ for the clause $\mathcal{C}_i$.

\begin{lemma} If $f = f(x_1,x_2,x_3)$ is a tropical quadratic polynomial
with 
\begin{equation} \label{eq:contained}
  \{ x \in \{0,1\}^3 \, : \, x_1 = 1 \text{ or } x_2 = 1 \text{ or }
       x_3 = 1 \} \subset \mathcal{T}(f) \, ,
\end{equation}
then $\mathcal{T}(f)$ also contains $(0,0,0)$.
\end{lemma}

\begin{proof}
Assume that there exists a tropical quadratic polynomial $f$
satisfying~\eqref{eq:contained} such that the minimum of the linear
forms at $(0,0,0)$ is attained only once. Let $l$ be the linear form
where the minimum is attained. Since $f$ is quadratic, $l$ depends
on at most two variables, and the exponents
of these variables are $1$ or $2$. Let $x_k$ be a variable which does
not occur in $l$, and let $x'$ be obtained from $x$
by switching $x_k$ from 0 to 1. Then the value of each linear form
at $x'$ is larger than or equal to the value of that linear form
of $x$. Since the value of the linear form $l$ at $x'$ is equal to
the value of $l$ at $x$, the minimum of all linear forms at $x'$ is
the same one as at $x$, and it is attained only once.
\end{proof}

In order to encode a clause with three positive literals into 
tropical quadrics, we embed it into higher-dimensional space
by introducing an additional variable $z$. 
Let $\mathcal{C}_i'$ be the Boolean formula
\[
  \mathcal{C}'_i \ = \ (y_{i_1} \vee y_{i_2} \vee \overline{z})
                      \wedge (y_{i_3} \vee z) 
                      \wedge (\overline{y_{i_3}} \vee \overline{z})
\]
in the variables $y_{i_1}, y_{i_2}, y_{i_3}, z$.
The last two clauses of this formula imply that any satisfying
assignment of $\mathcal{C}_i'$ has the property 
$y_{i_3} = \overline{z}$. Hence, there exists a satisfying
assignment for the original clause $\mathcal{C}_i$ if and only if
the formula $\mathcal{C}'_i$ can be satisfied. $\mathcal{C}_i'$
consists of one 3-clause that belongs to the case $p=2$ and
of two 2-clauses, which can be encoded into tropical geometry
as described above. Hence, there exist three tropical quadratic
polynomials $g_1,g_2,g_3$ in $y_{i_1}, y_{i_2}, y_{i_3}, z$ 
such that $\mathcal{C}_i$ can be satisfied if and only if
$\mathcal{T}(g_1)$, $\mathcal{T}(g_2)$, and $\mathcal{T}(g_3)$ have
a common point in $\{0,1\}^4$.

\medskip

For $p \in \{0,\ldots,3\}$ let $\#_p(\mathcal{C})$ denote the number of
clauses in the 3-\textsc{Sat} formula $\mathcal{C}$ with $p$
positive terms. Then the construction for the clauses yields
$k' := k + 2 \#_3(\mathcal{C})$ tropical hypersurfaces, which we denote
by $f_1, \ldots, f_{k'}$. Moreover, due to the additional auxiliary 
variables the actual number of total variables is 
$n' := n + \#_3(\mathcal{C})$. Let $P$ be the tropical prevariety
\begin{equation}
\label{eq:totaltropprevar}
  P \ = \ \mathcal{T}(h_1) \cap \ldots \cap \mathcal{T}(h_{n'}) 
  \cap \mathcal{T}(f_1) \cap \ldots \cap \mathcal{T}(f_{k'}) 
  \subset \R^{n'} \, .
\end{equation}

\begin{lemma} \label{le:canbesatisfied}
$P$ is nonempty if and only if $\mathcal{C}$ can be 
satisfied.
\end{lemma}

\begin{proof} Let 
$y \in \{\textsc{True},\textsc{False}\}^n$ 
be a satisfying assignment for
$\mathcal{C}$ and
$x := s(y) \in \{0,1\}^{n'}.$ 
By construction, $x$ is contained in all the structural hypersurfaces
and in all the clause surfaces.

Conversely, let $x \in P$. Since 
$x \in \bigcap_{i=1}^{n'} \mathcal{T}(h_i)$ we have $x \in \{0,1\}^{n'}$.
Set $y = s^{-1}(x) \in \{\textsc{True}, \textsc{False}\}^{n'}$.
Since $x$ is contained in all clause hypersurfaces representing the
clause $\mathcal{C}_i$,
the truth assignment $y$
satisfies the clause $\mathcal{C}_i$, $1 \le i \le k$.
Hence, $\mathcal{C}$ can be satisfied.
\end{proof}

All the polynomials in the construction of the tropical prevariety $P$
are of degree at most 2. Moreover, $P$ is a finite set and therefore 
even a tropical variety.
Since the reduction from 3-\textsc{Sat} to \textsc{Tropical Intersection}
is doable in polynomial time, this finishes the proof of
Lemma~\ref{le:tropintersection} and hence of the $\NP$-hardness statement
for \textsc{Tropical Intersection} in Theorem~\ref{th:hardness1}.

\begin{cor}
  $\#$\textsc{Connected Components} is $\#\mathcal{P}$-hard.
\end{cor}

\begin{proof}
It suffices to observe that the reduction given above is parsimonious,
i.e., the number of solutions of the tropical prevariety
is the number of satisfying assignments of the 3-$\textsc{Sat}$ formula
$\mathcal{C}$.
Since counting the number of satisfying assignments of a 3-$\textsc{Sat}$ 
formula is a ${\#}\classP$-hard problem \cite{valiant-78}, the statement follows.
\end{proof}

\begin{lemma} \label{le:tropintersectionnp}
\textsc{Tropical Intersection} $\in$ $\NP$.
\end{lemma}

\begin{proof}
We have to show that for every \textsc{Yes} instance of \textsc{Tropical
Intersection}, there exists a certificate of polynomial size, as
well as a polynomial time verification procedure for
these certificates.

Let $f_i$ be of the form $f_i = f_i(x_1, \ldots, x_n) = 
\bigoplus_{\alpha \in \mathcal{A}_i} c_\alpha \cdot x^{\alpha}$
for some support set $\mathcal{A}_i$, $1 \le i \le m$.
If there exists a point $z$ in the intersection of tropical varieties, 
then there exist 
$(\beta_1, \gamma_1) \in \mathcal{A}_1^2
, \ldots, (\beta_m, \gamma_m) \in \mathcal{A}_m^2$
such that the minimum in $f_i$ at $z$ is attained at the terms given by
$(\beta_i,\gamma_i)$. Hence,
$z$ is a solution of system of linear equations and 
inequalities
\[
  c_{\beta_i} + \sum_{j=1}^n \beta_{ij} x_j \ = \ 
  c_{\gamma_i} + \sum_{j=1}^n \gamma_{ij} x_j 
  \ \le \ c_{\alpha} + \sum_{j=1}^n \alpha_j x_j
  \quad \text{ for all }
  \alpha \in \mathcal{A}_i, \quad 1 \le i \le m \, .
\]
The size of this linear program 
is linear in the size of the input. Moreover, checking whether a given point $z$
is contained in a given tropical hypersurface can be done in polynomial
time. Consequently, checking whether $z$ is contained in the intersection
of tropical hypersurfaces can be done in polynomial time.
\end{proof}

Hence, by Lemmas~\ref{le:tropintersection} and~\ref{le:tropintersectionnp},
\textsc{Tropical Intersection} is $\NP$-complete. In contrast to this,
the following theorem provides a positive complexity result and
yields a linear program\-ming-based algorithm.

\begin{lemma} \label{le:fixedm}
If the number $m$ of tropical hypersurfaces
is a fixed constant, then
{\sc Tropical Intersection} can be solved in polynomial time.
\end{lemma}

\begin{proof}
Let $\mathcal{A}_i = \{ \alpha \in \N_0^n \, : \, 
  \sum_{j=1}^{d_j} \alpha_j \le d_i \}$
and $f_i = f_i(x_1, \ldots, x_n) = 
\bigoplus_{\alpha \in \mathcal{A}_i} c_\alpha \cdot x^{\alpha}$,
$1 \le i \le m$.
If $L$ denotes the binary encoding
length of the \textsc{Tropical Intersection}-instance then
the size $|\mathcal{A}_i|$ of $\mathcal{A}_i$ satisfies
$|\mathcal{A}_i| \le L$. Hence, for any $i \in \{1, \ldots, m\}$
the polynomial $f$ has at most $L$ terms, and thus there are at most
$\binom{L}{2}$ choices of two terms where the minimum in $f_i$ is attained.
Since there are at most $\binom{L}{2}^m$
choices of two terms in all the polynomials $f_1, \ldots, f_m$,
it suffices to show that for any fixed choice of two vectors 
$\beta_i, \gamma_i \in \mathcal{A}_i$ where the
minimum is attained in $f_i$, $1 \le i \le m$, the resulting linear 
program
\[
  c_{\beta_i} + \sum_{j=1}^n \beta_{ij} x_j \ = \ 
  c_{\gamma_i} + \sum_{j=1}^n \gamma_{ij} x_j 
  \ \le \ c_{\alpha} + \sum_{j=1}^n \alpha_j x_j
  \quad \text{ for all }
  \alpha \in \mathcal{A}_i, \quad 1 \le i \le m 
\]
can be solved in polynomial time. However, since the size of the
linear program is polynomial in the size of the input of
{\sc Tropical Intersection}, this follows from the polynomial
solvability of linear programming \cite{khachiyan-80}.
\end{proof}

\begin{lemma} \label{le:consistency}
{\sc Tropical Consistency} is co-$\NP$-hard.
This hardness persists if all polynomials are restricted to
be of degree at most~2.
\end{lemma}

\begin{proof} Since the empty set is a tropical variety,
it suffices to provide a polynomial time reduction from 
3-\textsc{Sat} with the following properties.
For every {\sc No}-instance of 3-\textsc{Sat} the constructed
tropical prevariety is the empty set. For every \textsc{Yes}-instance
of 3-\textsc{Sat}, the constructed prevariety is not a tropical
variety. In order to simplify notation, we assume from now on
that all clauses contain at most 2 positive literals, since otherwise
we can apply the same auxiliary construction as in the proof
of Lemma~\ref{le:tropintersection}.

We embed the construction from the proof of 
Lemma~\ref{le:canbesatisfied} into $\R^{n+1}$
by considering all polynomials formally to be polynomials 
in $n+1$ variables.
Since the definition of the structural hypersurfaces 
in~\eqref{eq:structuralclause1}
and the definition of the clause hypersurfaces 
do not depend on $x_{n+1}$, the embedding of the 
tropical prevariety $P$ from \eqref{eq:totaltropprevar}
into $\R^{n+1}$ gives a prevariety
$P' = P \times \R^n \subset \R$. 
Recall that the structural hypersurfaces are given by the 
tropical polynomials
$h_i = 0 \cdot x_i^2 \oplus 0 \cdot x_i \oplus 1$.
Let 
\[
  g_{i} \ = \
    0 \cdot x_i^2 \oplus 0 \cdot x_i \oplus 1 \oplus 1 \cdot  x_{n+1}
\]
for $1 \le i \le n$.
For $i=1$ and $n+1 = 2$, the tropical variety $\mathcal{T}(g_i)$ is shown in 
Figure~\ref{fi:embeddedreduction}.

\ifpictures
\begin{figure}[ht]
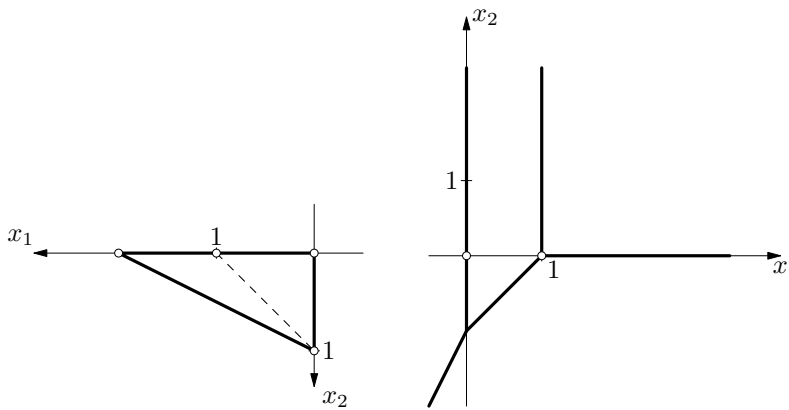

\[
  \begin{array}{c@{\qquad}c}
    \includegraphics[scale=1]{pictures/pictropcomp.4} &
    \includegraphics[scale=1]{pictures/pictropcomp.5} 
  \end{array}
\]

\caption{Newton polygon of the polynomial $g_{i}$
  for $i=1$, $n+1=2$ and $\mathcal{T}(g_{i})$.}
\label{fi:embeddedreduction}
\end{figure}
\fi

The intersection of $\mathcal{T}(g_i)$ and $\mathcal{T}(h_i)$ is
\[
  \mathcal{T}(g_{i}) \cap \mathcal{T}(h_i) \ = \ \left\{ 
  x \in \R^{n+1} \, : \, (x_i = 0 \text{ and } x_{n+1} \ge -1) \text{ or }
  (x_i = 1 \text{ and } x_{n+1} \ge 0) \right\} \, .
\]

If we imagine the $x_{n+1}$-axis to be pointing upwards, the intersection
of all structural hypersurfaces with the hypersurfaces
defined by $g_1, \ldots, g_n$ is a union of $2^n$ half
rays which are unbounded in the upward pointing directions,
\[
  \bigcap_{i=1}^n g_i \cap \bigcap_{i=1}^n h_i \ = \ 
  \left\{ \left( \{0,1\}^n \setminus (0, \ldots, 0) \right) \times \R_+ \, 
  \cup \, (0, \ldots, 0) \times \{x_{n+1} \in \R \, : \, x_{n+1} \ge -1\} \right\} \, .
\]

Using the same clause hypersurfaces as in the proof of
Lemma~\ref{le:tropintersection}, embedded into $\R^{n+1}$, we obtain the 
empty set for every {\sc No}-instance of 3-\textsc{Sat}.
Moreover, every \textsc{Yes}-instance of 3-\textsc{Sat}
gives a finite union of disjoint half rays which is not a tropical variety.
All the polynomials in the construction are of degree at most~2.
Since the reduction is polynomial time, the statement follows.
\end{proof}

\subsection{Connectivity\label{se:connectivity}}

In order to concentrate on the aspect of connectivity
(rather than a non-emptiness test in disguise), note that in the
definition of \textsc{Tropical Connectivity} 
we have excluded inputs leading to an empty prevariety.

\begin{lemma}
{\sc Tropical Connectivity} is co-$\NP$-hard.
This statement persists if the instances are restricted to
those where $\bigcap_{i=1}^m \mathcal{T}(f_i)$ is a tropical variety.
Moreover, this statement persists if all polynomials
are restricted to be of degree at most 3.
\end{lemma}

\begin{proof}
We choose a point $q$ which is always contained in the tropical variety
and modify the construction from the proof of Lemma~\ref{le:tropintersection}.
In order to achieve
that our choice of $q$ does not interfere with the remaining
construction, we embed the construction into $\R^{n+1}$, similar to 
the proof of Lemma~\ref{le:consistency}.

In the modification, the structural hypersurfaces are now given by
the polynomials
\[
  h_i'(x) \ = \ (t^0 \cdot x_i + t^1) \cdot (t^0 \cdot x_i + t^0) \cdot
            (t^0 \cdot x_i + t^2) \, , \quad 1 \le i \le n \, .
\]
Hence, the intersection of all hypersurfaces gives $\{0,1,2\}^{n+1}$.
By constructing additional polynomials of the form 
\[
    (t^0 \cdot x_i + t^1) \cdot (t^0 \cdot x_i + t^0) \cdot
            g_i'(x)
\]
with linear forms $g_i'$ as well as of the form
\[
    (t^0 \cdot x_{n+1} + t^0) \cdot (t^0 \cdot x_{n+1} + t^2)  \,,
\]
we can easily achieve that the intersection of
all theses hypersurfaces is the set
$\{0,1\}^n \times \{0\} \cup \{(2,2,2)\}$.

By multiplying all polynomials of the clause surfaces by
the polynomials
$t^0 \cdot x_{n+1} + t^2$, we can achieve that 
the point $(2,2,2)$ remains contained in the prevariety.
Note that the degree of all polynomials is increased by only 1.
Altogether, the constructed tropical prevariety $P$ is always
nonempty.
If the 3-$\textsc{Sat}$-formula can be satisfied,
then there are at least two connected components in $P$.
If the 3-$\textsc{Sat}$-formula cannot be satisfied then
$P$ has exactly one component.

All the constructed polynomials are of degree at most 3. The
resulting tropical prevariety is a finite set and therefore a
tropical variety. Moreover, the reduction is polynomial time.
\end{proof}

\begin{lemma} \label{le:fixedmconnectivity}
If the number $m$ of tropical prevarieties 
is a fixed constant, then
{\sc Tropical Connectivity} can be solved in polynomial time.
\end{lemma}

\begin{proof}
Similar to the proof of Lemma~\ref{le:fixedm}, for fixed $m$
we can compute
in polynomial time faces $F_1, \ldots, F_t$ of the polyhedral complex
$P = \bigcap_{i=1}^m \mathcal{T}(f_i)$ such that
$P = \bigcup_{i=1}^t F_i$. 
Hence, we can construct a graph
$G$ with vertices $F_1, \ldots, F_t$ in which two faces $F_i$ and $F_j$
are connected by an edge if and only if they
intersect. Then one computes the number of connected components
of $G$.
This can be done in polynomial time.
\end{proof}

The hardness result \ref{le:fixedmconnectivity} is
also contrasted by the statement that
linear tropical prevarieties are always connected.

\begin{lemma} \label{le:linearconnected}
Every nonempty linear tropical prevariety $P \subset \R^n$ is connected.
\end{lemma}

\begin{proof}
Let $P:= \bigcap_{i=1}^m \mathcal{T}(f_i)$ and $x,y \in P$. 
The notions $\oplus$, $\odot$, previously defined for scalars, 
can also be defined for vectors, by applying the operations componentwise.
With this notation 
it suffices to show that for every $\lambda, \mu \in \R$, the point 
$z := \lambda \odot x \oplus \mu \odot y$ is contained in 
each $\mathcal{T}(f_i)$, $1 \le i \le m$.
Fix an $i \in \{1, \ldots, m\}$, and let
$f_i = a_0 \oplus \bigoplus_{i=1}^n a_i \cdot x_i$. For convenience
of notation set $x_0 = y_0 = z_0 = 0$, 
and let $r$ be an index which minimizes 
$\{a_j + z_j \, : \, 0 \le j \le n\}$.
By definition of $z$, we have $z_r = \lambda + x_r$ or
$z_r = \mu + y_r$. Without loss of generality
we can assume $z_r = \lambda + x_r$. Note that then $a_r + x_r \le a_s + x_s$
for every index $s$.

Since $x \in \mathcal{T}(f_i)$,
there exists an index $s \neq r$ with $a_r + x_r = a_s + x_s$.
The definition of $z$ implies 
$a_s + z_s \le a_s + \lambda + x_s = a_r + \lambda + x_r = a_r + z_r$.
Hence, by the choice of $r$, $a_r + z_r = a_s + z_s$. 
In other words, the
minimum in $f_i$ is attained at least twice at the point~$z$,
i.e., $z \in \mathcal{T}(f_i)$.
 \end{proof}

\begin{rem}
Using the framework of tropical convexity 
from~\cite{develin-sturmfels-2004},
Lemma~\ref{le:linearconnected} also follows from the fact that
tropical hyperplanes are tropically convex 
\cite[Proposition~6]{develin-sturmfels-2004} in connection with
the observations that the intersection of tropically convex sets
is tropically convex and that tropically convex sets are connected.
\end{rem}

Statements~\ref{le:mconsistency}--\ref{le:linearconnected} prove
all claims in 
Theorems~\ref{th:hardness1}--\ref{th:sharpp}

\section{Related aspects on amoebas}

Our work is related to (and was partially inspired by) questions on
algorithmic complexity of basic problems on the
amoebas that were introduced in by Gel'fand, Kapranov,
Zelevinsky \cite{gkz-94}.
Let $I$ be an ideal in the ring $\R[x_1, \ldots, x_n]$ of
Laurent polynomials. Then the \emph{amoeba} of $I$ is defined by
the image of the complex subvariety $V(I) \subset (\C^*)$ 
under the mapping
\[
  \begin{array}{rcl}
    \text{Log}: (\C^*)^n & \to & \R^n \, , \\
                    z & \mapsto & (\log |z_1|, \ldots, \log |z_n|) \, , 
  \end{array}
\]
where $|\cdot|$ denotes the absolute value of a complex number
and $\log$ is the natural logarithm.
Since any hypersurface amoeba
contains a tropical variety (the so-called \emph{spine}) that is a
strong deformational retract of the amoeba (see, e.g.,
\cite[Theorem 2.6]{mikhalkin-survey-2004}),
algorithmic questions
on amoebas are closely related to algorithmic questions on
tropical varieties (see \cite[Chapter 9]{sturmfels-b2002}).

A central question by Einsiedler and Lind asks for an efficient algorithm
to test whether the complex amoeba of an ideal contains the
origin \cite{einsiedler-lind-question}. This comes
from applications in dynamical systems, where this test
determines whether a dynamical system has a
finiteness condition called expansiveness.
Only little is known about the computational hardness of
algorithmic problems on amoebas,
and the computational complexity of the membership problem
for amoebas (with rational input data
and the dimension being part of the input) is still open.
If the polynomials are given in sparse encoding (i.e.,
only the non-vanishing coefficients are listed in the input),
then the problem becomes $\NP$-hard even for $n=1$
(\cite{plaisted-84}, see also \cite{rojas-stella-2004}).
Recently, Rojas and Stella \cite{rojas-stella-2004} have established an algorithmic 
fewnomial theory providing further
hardness results for amoebas in sparse encoding (e.g., $\NP$-hardness
of deciding whether an amoeba intersects a coordinate hyperplane.)
For some Nullstellensatz-type algorithmic results see
\cite{purbhoo-2004}.

\bigskip

\noindent
{\bf Acknowledgment.} Thanks to Maurice Rojas for pointing out 
reference \cite{plaisted-84} and to the anonymous referees for very 
detailed criticism and comments.

\providecommand{\bysame}{\leavevmode\hbox to3em{\hrulefill}\thinspace}
\providecommand{\MR}{\relax\ifhmode\unskip\space\fi MR }
\providecommand{\MRhref}[2]{%
  \href{http://www.ams.org/mathscinet-getitem?mr=#1}{#2}
}
\providecommand{\href}[2]{#2}


\end{document}